\documentclass[a4paper,10pt]{amsart}
 \usepackage{amsfonts,mathrsfs,amsmath}
 \usepackage{amsrefs}
 \usepackage{amsthm}
 \usepackage{amssymb}
 \usepackage{enumerate}
 \usepackage{graphicx}
 \usepackage{tikz-cd}
 \usepackage{cleveref}
 \usepackage{cancel}
 \usepackage{float}
\allowdisplaybreaks

\theoremstyle{definition}

\theoremstyle{plain}

\newtheorem*{thm*}{Theorem}

\newtheorem*{cor*}{Corollary}

\newtheorem*{con*}{Conjecture}
\newtheorem*{frag*}{Question}
\newtheorem*{verm*}{Vermutung}

\newtheorem{theorem}{Theorem}
\newtheorem{proposition}{Proposition}

\newtheorem{corollary}{Corollary}
\newtheorem{lemma}{Lemma}
\theoremstyle{definition}

\usepackage{url}
\usepackage{kbordermatrix}
\theoremstyle{definition}

\usepackage{enumitem}
\numberwithin{equation}{section}

\makeatletter
\newcommand{\mylabel}[2]{#2\def\@currentlabel{#2}\label{#1}}
\makeatother

\usepackage[section]{algorithm}
\usepackage{algpseudocode}

\usepackage[a4paper,margin=2.5cm]{geometry}
\newcommand{\dd}{{\rm d}}
\newcommand{\mj}{\mathbb{J}}

\newcommand{\me}{\mathbb {E}}
\newcommand{\x}{\underline{x}}
\newcommand{\uu}{\underline{u}}
\newcommand{\uv}{\underline{v}}

\newcommand{\uz}{\underline{z}}

\newcommand{\RR}{\mathbb R}
\newcommand{\NN}{\mathbb N}

\newcommand{\calB}{\mathcal B}

\newcommand{\calD}{\mathcal D}

\newcommand{\calL}{\mathcal L}

\renewcommand{\kbldelim}{(}
\renewcommand{\kbrdelim}{)}
\renewcommand{\kbrowstyle}{\displaystyle}
\renewcommand{\kbcolstyle}{\displaystyle}

\DeclareMathOperator{\ext}{ext}

\DeclareMathOperator{\vol}{vol}

\DeclareMathOperator{\Card}{card}
\DeclareMathOperator{\opt}{\ast}
\DeclareMathOperator{\vertt}{Vert}

\newcommand{\cB}{{\mathcal B}}

\newcommand{\cD}{{\mathcal D}}

\usepackage{mathrsfs}
\usepackage{bbm}
\newcommand{\Int}{\operatorname{int}}

\title[Maximal signed volume for (multivariate) supermodular quasi-copulas]{Maximal signed volume for (multivariate) supermodular quasi-copulas}
\author{Matja\v z Omladi\v c${}^{1}$}
\address{Matja\v z Omladi\v c,
Faculty of Mathematics and Physics, University of Ljubljana  \& Institute of Mathematics, Physics and Mechanics, Ljubljana, Slovenia.} 
\email{matjaz@omladic.net}\thanks{${}^1$Supported by the ARIS (Slovenian Research and Innovation Agency)
research core funding No.\ P1-0448.}
\author{Martin Vuk}
\address{Martin Vuk, 
Faculty of Computer and Information Science, University of Ljubljana.}
\email{martin.vuk@fri.uni-lj.si}
\author{Alja\v z Zalar${}^{3}$}
\address{Alja\v z Zalar, 
Faculty of Computer and Information Science, University of Ljubljana  \& 
Faculty of Mathematics and Physics, University of Ljubljana  \&
Institute of Mathematics, Physics and Mechanics, Ljubljana, Slovenia.}
\email{aljaz.zalar@fri.uni-lj.si}
\thanks{${}^3$Supported by the ARIS (Slovenian Research and Innovation Agency)
research core funding No.\ P1-0288 and grants J1-50002, J1-60011.}

\newcommand{\comment}[1]{}
\newcommand{\mi}{\mathbb{I}}

\newcommand{\sign}{\mathrm{sign}}
\newcommand{\Vertices}{\operatorname{Vert}}

\begin{document}


\date{\today}

\begin{abstract} 
Copulas are the main tool of dependence modeling in statistics and quasi-copulas are their necessary companions. The latter appear, say, as infima or suprema of sets of copulas, they form a huge class and have some unpleasant properties. Their statistical interpretation is challenged by the fact that they may lead to negative volumes of some boxes. So, numerous applications call for an intermediate class and supermodular quasi-copulas is one of them having many useful properties. An excellent measure, Average Rectangular Volume (ARV in short), to clarify and position this class was proposed in the seminal paper by Anzilli and Durante, The average rectangular volume induced by supermodular aggregation functions, \emph{J.\ Math.\ Anal.\ Appl.} \textbf{555} (2026) 21 pp.
While supermodularity is a bivariate notion, its extension to $d$-variate case for $d>2$ was recently provided in a key paper by Arias-Garc\'{ı}a, Mesiar, and De Baets, The unwalked path between quasi-copulas and copulas: Stepping stones in higher dimensions, \emph{Int.\ J.\ of Appr.\ Reasoning}, \textbf{80} (2017) pp.\ 89–99. Here an alternative method to ARV is presented extendable to the multivariate case 
based on Maximal (in absolute value) Negative Volumes (MNV in short) on boxes thus helping practitioners when seeking the right (quasi-)copula for their problem. Observe that these volumes on copulas are zero, while their values on quasi-copulas, depending on $d$, has been a long-open problem solved only recently.  We present a nontrivial extension of this solution to exhibit the main goal of this paper, a measure that clarifies and 
positions the classes considered based on MNV.

\end{abstract}

\maketitle

\section{Introduction}

Copula \cite{durante2015principles} is a joint distribution with uniform margins, but when Sklar \cite{Sklar1959229} discovered in 1959 that we could get any distribution on earth when we inserted appropriate univariate distributions as margins into a copula, it became a universal 
device to solve dependence problems. There are many applications, though, where quasi-copulas are also needed, they form a pretty big class and have some unpleasant properties. They appear, say, as infima or suprema of sets of copulas, but their statistical interpretation is challenged by the fact that they may lead to negative volumes of some boxes. So, numerous applications call for an intermediate class and supermodular quasi-copulas is one of them having many useful properties. 

An excellent method, Average Rectangular Volume (ARV in short), to aggregate the knowledge contained in a quasi-copula of this class was proposed by Anzilli and Durante in their seminal paper \cite{AD26} in 2026. 
While supermodularity is a bivariate notion \cite{DMPS07}, its extension to $d$-variate case for $d>2$ was recently provided in a key paper by Arias-Garc\'{ı}a et al.\ \cite{AGMB17} in 2017 and we want to find a parallel notion to ARV for
these classes between the classes of quasi-copulas they are contained in and copulas they contain. There seems to be no obvious extension of the method presented in \cite{AD26}, so, we develop a novel method based on maximal (in absolute value) negative volume on boxes. Observe that this value on copulas is zero, while its value on quasi-copulas, depending on $d$, has been a long-open problem solved only recently. Namely, the seminal paper \cite{ARIASGARCIA20201} by Arias-Garc\'{ı}a et al.\  
which is in some sense a 2020 sequel of \cite{AGMB17} presents a survey of the most significant results and open problems in quasi-copula theory. Their Open Problem \#5 asks for the maximal negative and maximal positive mass \emph{over all boxes and over all quasi-copulas}, for the characterization of the boxes attaining these extremes, and for characterization of quasi-copulas with this property. Here, \textit{maximal negative} refers to the negative volume with the largest absolute magnitude.
The problem seems to have been first proposed in 2002 by Nelsen et al.\ \cite{Nelsen2002}, the trivariate solution was first presented in 2007 by De Beats et al.\ \cite{BMUF07} and the problem was solved in papers \cite{UF23} in 2023 and \cite{BeOmVuZa,OmVuZa25} in 2025/26. 

The notion of \textit{supermodularity} (see, e.g., \cite[Definition~2.1]{KMM11}) coincides with the 2-increasing property in the bivariate case. 
In dimensions $d\geq\emph{} 2$, beyond the class of 2-increasing functions, there exists a hierarchy of intermediate classes satisfying the $k$-increasing property for $k=2,\ldots,d$. These classes were introduced in \cite[Section~5]{AGMB17} (see also \cite[Section~6.3]{ARIASGARCIA20201}), extending the volume-based characterization of quasi-copulas to $k$-dimensional slices of the $d$-dimensional box. In this way one obtains special subclasses of quasi-copulas: they are strictly contained within the class of all quasi-copulas (corresponding to the 1-increasing property) and contain all copulas (corresponding to the $d$-increasing property).
Consequently, the notion of $k$-increasingness can be viewed as a natural generalization of supermodularity to higher dimensions as we mentioned earlier. 
The Maximal (in absolute value) Negative Volumes (MNV in short) on boxes and their counterpart Maximal Positive Volumes (MPV in short) bring us to their tie Maximal Signed Volumes (MSV for short), thus providing a possible parallel of ARV.

The manuscript is organized as follows. In Section~\ref{sec:main-theorems} we present our main results (Theorems~\ref{sol:general-increasing} and \ref{sol:general-increasing-maximal}), together with a numerical analysis of the extremal values (Tables~\ref{table:minimal-volume} and \ref{table:maximal-volume}). Section~\ref{preparatory} formulates the maximal-volume problems as linear programs (Proposition~\ref{prel:k-prop}), simplifies them using a symmetrization argument (Corollary~\ref{cor:sim-k}), introduces new sets of variables (Proposition~\ref{simplify-further-general}), and derives the dual linear programs (Proposition~\ref{simplify-further-general-v2}). In Section~\ref{sec:proofs} we prove Theorems~\ref{sol:general-increasing} and \ref{sol:general-increasing-maximal} by solving these dual programs exactly.
Finally, in Section \ref{sec:concluding} we summarize the paper and give some directions for further work.

\section{Statements of the main theorems}
\label{sec:main-theorems}

{In this section we introduce the notation and state our main results, i.e., Theorem \ref{sol:general-increasing} solves the maximal negative volume problem for $k$-increasing quasi-copulas, while Theorem \ref{sol:general-increasing-maximal} solves the maximal positive volume problem. Concrete numerical solution examples are also presented in the form of tables. (See Tables \ref{table:minimal-volume} and \ref{table:maximal-volume} for the minimal and the maximal volume question, respectively.)}\\

Let $\calD\subseteq [0,1]^d$ be a non-empty set and $d\in \NN$, $d\geq 2$.
We say that a function 
    $F:\cD\to [0,1]$:
\begin{itemize}
    \item is \emph{$d$-increasing} if for any $d$-box
$
\cB:=\prod_{i=1}^d [a_i,b_i]\subseteq \cD
$
it holds that
\[
V_F(\cB):=\sum_{z\in \Vertices(\cB)} (-1)^{S(z)}\,F(z)\ \ge 0,
\]
where $\operatorname{Vert}(\cB)$ stands for the vertices of $\cB$ and $S(z)=\Card\{\,j\in\{1,2,\dots,d\}\colon z_j=a_j\,\}$. The quantity $V_F(\cB)$ is called the \emph{$F$-volume} of $\cB$. Here $\Card A$ denotes the cardinality of the set $A$.
    \item
    satisfies the \emph{boundary condition}
    if for $\underline u:=(u_1, \ldots, u_d)\in \calD$ the following hold:
        \begin{align*} 
        &(a)\quad 
        \text{If }        
        \underline 
u=(u_1, \ldots, u_{i-1}, 0, u_{i+1}, \ldots, u_d)
        \text{ for some }i, 
        \text{then }F(\underline u)=0.\\
        &(b)\quad 
        \text{If }
        \underline u=(1, \ldots, 1, u_i, 1, \ldots, 1)
        \text{ for some }i, 
        \text{then }F(\underline u)=u_i.
        \end{align*}
    \item
    satisfies the \emph{monotonicity condition} 
    if 
    it is nondecreasing in every variable, i.e., for each $i=1,\ldots,d$ and each pair of $d$-tuples
    \begin{align*}
    \begin{split}
    \uu&:=(u_1, \ldots,u_{i-1},u_i,u_{i+1},\ldots,u_d)\in \calD,\\
    \widetilde\uu&:=(u_1, \ldots,u_{i-1},\widetilde u_i,u_{i+1},\ldots,u_d)\in \calD,
    \end{split}
    \end{align*}
    such that $u_i\leq \widetilde u_i$, it follows that
    $F(\uu)\leq F(\widetilde\uu)$.
    \item
    satisfies the \emph{Lipschitz condition}
    if for
        given $d$--tuples 
        $(u_1, \ldots, u_d)$ and $(v_1, \ldots, v_d)$ 
        in $\calD$, it holds that
        \begin{align*}
        |F(u_1, \ldots, u_d) - F(v_1, \ldots, v_d)| \leq \sum_{i = 1}^{d}|u_i - v_i|.
        \end{align*}
\end{itemize}
If $\calD=[0,1]^d$ and $F$ satisfies 
the boundary, the monotonicity and the Lipschitz conditions, then $F$ is called a \emph{$d$-variate quasi-copula}. We will omit the dimension $d$ when it is clear from the context and write quasi-copula for short.

Let $\cB:=\prod_{j=1}^d [x_j,y_j]$, where $0\leq x_j\leq y_j\leq 1$ for each $j$.
Next we recall the definition of a $k$-dimensional section of a function $F:\cB\to [0,1]$
(see \cite[Definition 5]{AGMB17}).
For any $\underline a:=(a_1,\ldots,a_d) \in \cB$ and a set of indices $A \subseteq \{1,2,\ldots,d\}$ with 
$0 < \Card A = k \leq d$, the
\emph{$k$-dimensional slice} (or \emph{$k$-slice}) $\cB_{\underline a,A}$ 
of $\cB$ with fixed values given by $\underline{a}$ in the positions determined by $A$,  is defined by 
$$
\cB_{\underline a,A}:=
\{(z_1,\ldots,z_d)\colon
z_j=a_j \; \text{if }j\notin A
\text{ and }z_j\in [x_j,y_j]\text{ if }j\in A.
\}
$$
In particular, if $\underline a$ is a vertex of $\cB$, then
$\cB_{\underline a,A}$ is called a \emph{$k$-dimensional face} (or \emph{$k$-face}) of $\cB$.
A \emph{$k$-dimensional section} (or \emph{$k$-section}) of $F$ with fixed values given by $\underline{a}$ in the positions determined by $A$, 
is the function
\begin{equation}
\label{def:k-section}
F_{\underline{a},A} : B_{\underline a,A} \to [0,1],
\end{equation}
defined by
\[
F_{\underline{a},A}(\underline{z}) := F(\underline{w}),
\quad \text{where} \quad
w_j =
\begin{cases}
z_j, & \text{if } j \in A,\\
a_j, & \text{if } j \notin A.
\end{cases}
\]
We say that a function $F:\cB \to [0,1]$ is \emph{$k$-dimensionally-increasing} (or \emph{$k$-increasing}), with $k \in \{1, \ldots, d\}$, if any of its $k$-sections is $k$-increasing.

We denote by $\mathrm{MNV}(d,k)$ the maximal negative volume $V_{Q}(\cB)$, 
i.e., the largest in absolute value among negative ones,
of some $d$-box $\cB$ over all $k$-increasing $d$-variate quasi-copulas $Q$.
The following result provides an explicit formula for $\mathrm{MNV}(d,k)$.


\begin{theorem}
\label{sol:general-increasing}
    Assume the notation above.
    Let $d,k\in \NN$ and $2\leq k\leq d$.
    Then
\begin{equation*}
\mathrm{MNV}(d,k)=\min\Big(0,\min_{i=k,\ldots,d}
        (-1)^{d-i}\binom{d-k}{i-k}\left(\gamma_i^{(k)}\right)^{-1}\Big),
\end{equation*}
where $\gamma_i^{(j)}$ are defined recursively by
\begin{align}
\label{recursiveness-v2}
\begin{split}
    \gamma_{i}^{(j)}&=
        \left\{
        \begin{array}{rl}
            d+1-i & \text{for }j=2\text{ and }i=2,\ldots,d,\\[0.2em]
            \gamma_i^{(j-1)}& \text{for }j= 3,\ldots,k\text{ and }i= 2,\ldots,j-2,\\[0.2em]
            \sum_{\ell=i}^{d}\gamma_{\ell}^{(j-1)}& 
                \text{for }j= 3,\ldots,k\text{ and }
                    i= j-1,\ldots,d.
        \end{array}
        \right.\\
\end{split}
\end{align}

Let $i_0\in \{k,\ldots,d\}$ be such that
$$\mathrm{MNV}(d,k)
=(-1)^{d-i_0}\binom{d-k}{i_0-k}\left(\gamma_{i_0}^{(k)}\right)^{-1}.$$ 
One of the realizations of $Q$ and $\cB$ such that $V_Q(\cB) = \mathrm{MNV}(d,k)$
has the following properties:
\begin{equation}
\label{one:realiziation}
        \cB=\left[a,1\right]^d,
        \qquad
        a=\frac{\alpha_{i_0}^{(k)}}{\gamma_{i_0}^{(k)}},
        \qquad
        Q(z) = q_{d-S(z)},\quad z\in\Vertices(\cB),
        \end{equation}
where $q_i$ are defined as follows
        \begin{equation}
        \label{one:recursive_qi}
        q_i=
        \left\{
        \begin{array}{rl}
            0,& \text{for }i= 0,\ldots,i_0-1,\\[0.3em]
            \left(\gamma_{i_0}^{(k)}\right)^{-1},&     \text{for }i=i_0,\\[0.3em]
            \sum_{j=1}^{k}\binom{k}{j}(-1)^{j+1}q_{i-j},& 
            \text{for }i= i_0+1,\ldots,d,
        \end{array}
        \right.
\end{equation}
and $\alpha_i^{(j)}$ are defined recursively by 
\begin{align}
\label{recursiveness-v3}
    \alpha_{i}^{(j)}
        &=
        \left\{
        \begin{array}{rl}
            1 & \text{for }j=2\text{ and }i= 2,\ldots,d,\\[0.2em]
            \alpha_i^{(j-1)}& \text{for }j= 3,\ldots,k\text{ and }i= 2,\ldots,j-2,\\[0.2em]
            \sum_{\ell=i}^{d}\alpha_{\ell}^{(j-1)}& 
                \text{for }
                j= 3,\ldots,k
                \text{ and }
                i= j-1,\ldots,d.
        \end{array}
        \right.
\end{align}
    This realization of $Q$ on 
    $\left\{a,1\right\}^d$ indeed extends to a 
    quasi-copula $Q:[0,1]^d\to [0,1]$
    by Proposition \ref{prel:k-prop} below.
\end{theorem}

We denote by $\mathrm{MPV}(d,k)$ the maximal positive volume $V_{Q}(\cB)$
of some $d$-box $\cB$ over all $k$-increasing $d$-variate quasi-copulas $Q$.
The following result provides an explicit formula for $\mathrm{MPV}(d,k)$.


\begin{theorem}
\label{sol:general-increasing-maximal}
    Assume the notation above.
    Let $d,k\in \NN$ and $2\leq k\leq d$.
    Then
\begin{equation*}
\mathrm{MPV}(d,k)
=\max_{i=k,\ldots,d}
        (-1)^{d-i}\binom{d-k}{i-k}(\gamma_i^{(k)})^{-1},
\end{equation*}
where $\gamma_i^{(j)}$ are defined recursively by
\eqref{recursiveness-v2}.

Let $i_0\in \{k,\ldots,d\}$ 
be such that
$$
\mathrm{MPV}(d,k)
=(-1)^{d-i_0}\binom{d-k}{i_0-k}(\gamma_{i_0}^{(k)})^{-1}.$$ 
One of the realizations of $\cB$ and $Q$ is 
as in \eqref{one:realiziation}.
\end{theorem}

We may also consider the maximal signed volume as the pair $\mathrm{MSV}(d,k)=(-\mathrm{MNV}(d,k),\mathrm{MPV}(d,k))$; observe that $\mathrm{MSV}(2,2)$ may be seen as a possible statistical counterpart to ARV of \cite{AD26}.\\[1mm]

In Table \ref{table:minimal-volume} and \ref{table:maximal-volume}\footnote{The numerical analysis in arithmetic over $\mathbb{Q}$ was performed using the software tool \textit{Mathematica} \cite{ram2024}. The source code is available at  \url{https://github.com/ZalarA/Quasi-copulas-k-increasing-extreme-volumes}.} we give 
the values of 
$\mathrm{MNV}(d,k)$ and 
$\mathrm{MPV}(d,k)$
for dimensions up $d$ to 15. The cases $k\geq 2$ are computed according to Theorems \ref{sol:general-increasing} and \ref{sol:general-increasing-maximal} above, while the case $k=1$ was studied in our previous work \cite[Theorems 1 and 2]{OmVuZa25}.

\begin{table}[h!]
      \caption{
      Explicit values of $\mathrm{MNV}(d,k)$
      for $2\leq d\leq 15$ and $1\leq k\leq d$.
      }
{
\renewcommand{\arraystretch}{1.4} 
\begin{tabular}{|c|c|c|c|c|c|c|c|c|c|c|c|c|c|c|}
\hline
$k\backslash d$ & 2 & 3 & 4 & 5 & 6 & 7 & 8 & 9 & 10 & 11 & 12 & 13 & 14 & 15 \\[4pt]
\hline
1 & $-\tfrac{1}{3}$ & $-\tfrac{4}{5}$ & $-\tfrac{9}{7}$ & $-\tfrac{32}{13}$ & $-\tfrac{75}{16}$ & $-\tfrac{19}{2}$ & $-\tfrac{55}{3}$ & $-37$ & $-\tfrac{209}{3}$ & $-\tfrac{251}{2}$ & $-\tfrac{791}{3}$ & $-\tfrac{1583}{3}$ & $-\tfrac{3002}{3}$ & $-\tfrac{7435}{4}$ \\[4pt]
\hline
2 & $0$ & $-\tfrac{1}{2}$ & $-1$ & $-\tfrac{3}{2}$ & $-2$ & $-\tfrac{5}{2}$ & $-5$ & $-\tfrac{35}{4}$ & $-14$ & $-21$ & $-42$ & $-77$ & $-132$ & $-\tfrac{429}{2}$ \\[4pt]
\hline
3 &  & $0$ & $-\tfrac{1}{3}$ & $-\tfrac{2}{3}$ & $-1$ & $-\tfrac{4}{3}$ & $-\tfrac{5}{3}$ & $-2$ & $-\tfrac{7}{2}$ & $-\tfrac{28}{5}$ & $-\tfrac{42}{5}$ & $-12$ & $-22$ & $-\tfrac{264}{7}$ \\[4pt]
\hline
4 &  &  & $0$ & $-\tfrac{1}{4}$ & $-\tfrac{1}{2}$ & $-\tfrac{3}{4}$ & $-1$ & $-\tfrac{5}{4}$ & $-\tfrac{3}{2}$ & $-\tfrac{7}{4}$ & $-\tfrac{14}{5}$ & $-\tfrac{21}{5}$ & $-6$ & $-\tfrac{33}{4}$ \\[4pt]
\hline
5 &  &  &  & $0$ & $-\tfrac{1}{5}$ & $-\tfrac{2}{5}$ & $-\tfrac{3}{5}$ & $-\tfrac{4}{5}$ & $-1$ & $-\tfrac{6}{5}$ & $-\tfrac{7}{5}$ & $-\tfrac{8}{5}$ & $-\tfrac{12}{5}$ & $-\tfrac{24}{7}$ \\[4pt]
\hline
6 &  &  &  &  & $0$ & $-\tfrac{1}{6}$ & $-\tfrac{1}{3}$ & $-\tfrac{1}{2}$ & $-\tfrac{2}{3}$ & $-\tfrac{5}{6}$ & $-1$ & $-\tfrac{7}{6}$ & $-\tfrac{4}{3}$ & $-\tfrac{3}{2}$ \\[4pt]
\hline
7 &  &  &  &  &  & $0$ & $-\tfrac{1}{7}$ & $-\tfrac{2}{7}$ & $-\tfrac{3}{7}$ & $-\tfrac{4}{7}$ & $-\tfrac{5}{7}$ & $-\tfrac{6}{7}$ & $-1$ & $-\tfrac{8}{7}$ \\[4pt]
\hline
8 &  &  &  &  &  &  & $0$ & $-\tfrac{1}{8}$ & $-\tfrac{1}{4}$ & $-\tfrac{3}{8}$ & $-\tfrac{1}{2}$ & $-\tfrac{5}{8}$ & $-\tfrac{3}{4}$ & $-\tfrac{7}{8}$ \\[4pt]
\hline
9 &  &  &  &  &  &  &  & $0$ & $-\tfrac{1}{9}$ & $-\tfrac{2}{9}$ & $-\tfrac{1}{3}$ & $-\tfrac{4}{9}$ & $-\tfrac{5}{9}$ & $-\tfrac{2}{3}$ \\[4pt]
\hline
10 &  &  &  &  &  &  &  &  & $0$ & $-\tfrac{1}{10}$ & $-\tfrac{1}{5}$ & $-\tfrac{3}{10}$ & $-\tfrac{2}{5}$ & $-\tfrac{1}{2}$ \\[4pt]
\hline
11 &  &  &  &  &  &  &  &  &  & $0$ & $-\tfrac{1}{11}$ & $-\tfrac{2}{11}$ & $-\tfrac{3}{11}$ & $-\tfrac{4}{11}$ \\[4pt]
\hline
12 &  &  &  &  &  &  &  &  &  &  & $0$ & $-\tfrac{1}{12}$ & $-\tfrac{1}{6}$ & $-\tfrac{1}{4}$ \\[4pt]
\hline
13 &  &  &  &  &  &  &  &  &  &  &  & $0$ & $-\tfrac{1}{13}$ & $-\tfrac{2}{13}$ \\[4pt]
\hline
14 &  &  &  &  &  &  &  &  &  &  &  &  & $0$ & $-\tfrac{1}{14}$ \\[4pt]
\hline
15 &  &  &  &  &  &  &  &  &  &  &  &  &  & $0$ \\[4pt]
\hline
\end{tabular}
}
\label{table:minimal-volume}
\end{table}
\vfill

\clearpage

\begin{table}[H]
      \caption{
      Explicit values of $\mathrm{MPV}(d,k)$
      for $2\leq d\leq 15$ and $1\leq k\leq d$.
      }
{
\renewcommand{\arraystretch}{1.4}
\begin{tabular}{|c|c|c|c|c|c|c|c|c|c|c|c|c|c|c|}
\hline
$k\backslash d$ & $2$ & $3$ & $4$ & $5$ & $6$ & $7$ & $8$ & $9$ & $10$ & $11$ & $12$ & $13$ & $14$ & $15$ \\
\hline
$1$ & $1$ & $1$ & $2$ & $\tfrac{7}{2}$ & $\tfrac{11}{2}$ & $\tfrac{31}{3}$ & $19$ & $\tfrac{71}{2}$ & $\tfrac{211}{3}$ & $\tfrac{421}{3}$ & $\tfrac{793}{3}$ & $\tfrac{1915}{4}$ & $\tfrac{3004}{3}$ & $\tfrac{6007}{3}$ \\
\hline
$2$ & $1$ & $1$ & $1$ & $1$ & $2$ & $\tfrac{10}{3}$ & $5$ & $7$ & $14$ & $\tfrac{126}{5}$ & $42$ & $66$ & $132$ & $\tfrac{1716}{7}$ \\
\hline
$3$ &  & $1$ & $1$ & $1$ & $1$ & $1$ & $\tfrac{5}{3}$ & $\tfrac{5}{2}$ & $\tfrac{7}{2}$ & $\tfrac{14}{3}$ & $\tfrac{42}{5}$ & $14$ & $22$ & $33$ \\
\hline
$4$ &  &  & $1$ & $1$ & $1$ & $1$ & $1$ & $1$ & $\tfrac{3}{2}$ & $\tfrac{21}{10}$ & $\tfrac{14}{5}$ & $\tfrac{18}{5}$ & $6$ & $\tfrac{66}{7}$ \\
\hline
$5$ &  &  &  & $1$ & $1$ & $1$ & $1$ & $1$ & $1$ & $1$ & $\tfrac{7}{5}$ & $\tfrac{28}{15}$ & $\tfrac{12}{5}$ & $3$ \\
\hline
$6$ &  &  &  &  & $1$ & $1$ & $1$ & $1$ & $1$ & $1$ & $1$ & $1$ & $\tfrac{4}{3}$ & $\tfrac{12}{7}$ \\
\hline
$7$ &  &  &  &  &  & $1$ & $1$ & $1$ & $1$ & $1$ & $1$ & $1$ & $1$ & $1$ \\
\hline
$8$ &  &  &  &  &  &  & $1$ & $1$ & $1$ & $1$ & $1$ & $1$ & $1$ & $1$ \\
\hline
$9$ &  &  &  &  &  &  &  & $1$ & $1$ & $1$ & $1$ & $1$ & $1$ & $1$ \\
\hline
$10$ &  &  &  &  &  &  &  &  & $1$ & $1$ & $1$ & $1$ & $1$ & $1$ \\
\hline
$11$ &  &  &  &  &  &  &  &  &  & $1$ & $1$ & $1$ & $1$ & $1$ \\
\hline
$12$ &  &  &  &  &  &  &  &  &  &  & $1$ & $1$ & $1$ & $1$ \\
\hline
$13$ &  &  &  &  &  &  &  &  &  &  &  & $1$ & $1$ & $1$ \\
\hline
$14$ &  &  &  &  &  &  &  &  &  &  &  &  & $1$ & $1$ \\
\hline
$15$ &  &  &  &  &  &  &  &  &  &  &  &  &  & $1$ \\
\hline
\end{tabular}
}
\label{table:maximal-volume}
\end{table}

\section{Towards the proof of Theorems \ref{sol:general-increasing} and \ref{sol:general-increasing-maximal}}
\label{preparatory}

In this section, we formulate the maximal volume problems as linear programs (see Proposition~\ref{prel:k-prop}). We then simplify these programs in two steps: first, by applying the symmetrization trick (see Corollary~\ref{cor:sim-k}), and second, by introducing new sets of variables (see Proposition~\ref{simplify-further-general}).
We note that the symmetrization trick played a crucial role in solving the maximal volume problems for quasi-copulas in our previous work~\cite{OmVuZa25}, whereas the introduction of new sets of variables is the main novelty of the present paper, as it substantially reduces the complexity of the resulting linear programs.
Finally, we derive the corresponding dual linear programs (see Proposition~\ref{simplify-further-general-v2}), which will be used in the next section to prove Theorems~\ref{sol:general-increasing} and \ref{sol:general-increasing-maximal}.
\\


Fix $d\in \NN$. 
For multi-indices 
$$
\mi=(\mi_1,\ldots,\mi_d)\in \{0,1\}^d\quad \text{and} 
\quad
\mj=(\mj,\ldots,\mj_d)\in \{0,1\}^d$$ 
let
    $$\mj-\mi=(\mj-\mi_1,\ldots,\mj_d-\mi_d)\in \{-1,0,1\}^d$$ 
stand for their usual coordinate-wise difference.
Let 
$\me^{(\ell)}$ stand for the multi-index with the only non-zero coordinate the $\ell$--th one, which is equal to 1.
\newcommand{\rel}{\prec}
For each $\ell=1,\ldots,d$ we define a relation on $\{0,1\}^d$ by
    $$
        \mi \rel_\ell \mj 
        \quad \Leftrightarrow \quad
        \mj-\mi=\me^{(\ell)}.
    $$
We write
    $$
        \mi \rel \mj 
        \quad \Leftrightarrow \quad
        \mi \rel_\ell \mj\;\; \text{ for some }\ell\in \{1,2,\ldots,d\}.
    $$

For a point $\x=(x_1,\ldots,x_d)\in \RR^d$ we define the functions
\begin{align*}
    G_d:\RR^d\to \RR,\quad 
    G_d(\x)
    &:=\sum_{i=1}^d x_i-d+1,\\
    H_d:\RR^d\to \RR,\quad
    H_d(\x)
    &:=\min\{x_1, x_2, \ldots, x_d\}.
\end{align*}

Let $Q$ be a quasi-copula and $\cB=\prod_{i=1}^d [a_i,b_i]\subseteq [0,1]^d$ a $d$-box with $a_i\leq b_i$
for each $i$.
We will use multi-indices
of the form $\mi := (\mi_1, \mi_2 \ldots \mi_d) \in \{0, 1\}^d$ to index 
$2^d$ vertices 
of $\cB$.
Let $\|\mi\|_1:=\sum_{j=1}^d \mi_j$ be the $1$-norm of $\mi$.
We write 
\begin{equation*}
    x_\mi:=((x_\mi)_1,\ldots,(x_\mi)_d) 
\end{equation*}
to denote the vertex with coordinates
\begin{equation}
\label{def:xI-v2}
    (x_\mi)_k = \begin{cases}
    a_k,\quad \text{if }\mi_k=0,\\
    b_k,\quad \text{if }\mi_k=1.
    \end{cases}
\end{equation}
Let us denote the value of $Q$ in the point $x_\mi$
by 
\begin{equation*}
    q_\mi := Q(x_\mi).
\end{equation*}
We write $\sign(\mi) := (-1)^{d-\|\mi\|_1}$. 
In the notation above, the {$Q$-volume} of $\cB$ is equal to
\begin{equation*}
    V_Q(\cB) = \sum_{\mi\in \{0, 1\}^d} \sign(\mi) \cdot q_\mi.
\end{equation*}
Let $F$ be a face of the unit $d$-cube $[0, 1]^d$. We denote with 
$$
\mi(F):=\{\mi\colon \mi\in\Vertices(F)\}\subseteq\{0, 1\}^d
$$
the set of all multi-indices, corresponding to the vertices of the face $F$. Let 
$$
\mi_{\max}(F)
:=(
    \max_{\mi \in \mi(F)} \mi_1,
    \max_{\mi \in \mi(F)} \mi_2,
    \ldots, 
    \max_{\mi \in \mi(F)} \mi_d
    )
$$ 
be the vertex of $F$ such that each of its coordinates is the largest among all vertices.
We define
$$\sign_{F}(\mi) := (-1)^{{\|\mi_{\max}(F)-\mi\|}_1}.$$

We will prove 
Theorems \ref{sol:general-increasing} 
and \ref{sol:general-increasing-maximal}
in several steps using linear programming as the main tool.

\begin{proposition}
\label{prel:k-prop}
For any $k\in \{1, 2,\ldots, d\}$ define the following linear program
\begin{align}
\label{LP-general}
\begin{split}
\min_{
\substack{
    a_1,\ldots,a_d,\\
    b_1,\ldots,b_d,\\
    q_{\mi} \text{ for }\mi\in \{0,1\}^d
}}
&\hspace{0.2cm} \sum_{\mi\in \{0,1\}^d} \sign(\mi) q_\mi,\\
\text{subject to }
&\hspace{0.2cm} 
    0\leq a_i< b_i\leq 1\quad \text{for all }i=1,2,\ldots,d,\\
&\hspace{0.2cm} 
    q_{\mj} - q_{\mi} \le b_\ell - a_\ell
\quad 
    \text{for all }\ell=1,2,\ldots,d
    \text{ and all }\mi\rel_\ell\mj,\\
&\hspace{0.2cm}
    0\leq \sum_{\mj\in \mi(F)}
    \sign_{F}(\mj)\cdot q_\mj
    \quad
    \text{for all $j$-dimensional faces $F$ of }[0,1]^d \\
&\hspace{4.5cm}    
    \text{ and for all } j= 1,\ldots,k,\\
&\hspace{0.2cm} 
    \max\{0,G_d(x_{\mi})\} 
    \le q_\mi \le 
    H_d(x_{\mi})
    \quad \text{for all }\mi\in \{0,1\}^d.
\end{split}
\end{align}
Let $\mi^{(1)},\ldots,\mi^{(2^d)}$ be some order of all multi-indices $\mi\in\{0,1\}^d$.
    If there exists an optimal solution 
    $$(a_1^\ast,\ldots,a_d^\ast,
    b_1^\ast,\ldots,b_d^\ast,
    q_{\mi^{(1)}},\ldots,q_{\mi^{(2^d)}}
    )$$
    to \eqref{LP-general}, which
    satisfies 
    \begin{equation}
        \label{bi-equal-to-1}
            b_1^\ast=\ldots=b_d^\ast=1,
    \end{equation}
    then the optimal value of \eqref{LP-general} is the maximal negative volume of some box over all $k$-increasing $d$-variate quasi-copulas.

    Moreover, if there exists an optimal solution 
    to \eqref{LP-general} where $\min$ is replaced with $\max$, which
    satisfies 
    \eqref{bi-equal-to-1},
    then the optimal value of \eqref{LP-general} is the maximal positive volume of some box over all $k$-increasing $d$-variate quasi-copulas.
\end{proposition}

In the proof of Proposition \ref{prel:k-prop} we will use the fact that the $k$-increasing property of a multilinear function on a $d$-box follows from the $k$-increasing property on all $k$-faces of the box (see Lemma \ref{prop:ml-lipschitz} below).

\begin{lemma}\label{prop:ml-lipschitz}
    Let $F$
    be a multilinear function on a $d$-box $\cB=\prod_{i=1}^d[x_i, y_i]\subseteq [0,1]^d$. 
    The following statements are equivalent:
\begin{enumerate}
    \item\label{prop:ml-lipschitz-pt1} 
        $F$ is $k$-increasing.
    \item\label{prop:ml-lipschitz-pt2} 
        $F$ is $k$-increasing on each $k$-face of $\cB$.
\end{enumerate} 
\end{lemma}

\begin{proof}
    The nontrivial implication is 
    $\eqref{prop:ml-lipschitz-pt2}\Rightarrow\eqref{prop:ml-lipschitz-pt1}$.
    We will use induction on $\ell$ 
    to show that $F$ is 
    $k$-increasing 
    on each $\ell$-dimensional face of $\cB$ for $\ell=k,k+1,\ldots,d$.
    For $\ell=d$ we get that $F$ is $k$-increasing on the whole box $\cB$, proving
    \eqref{prop:ml-lipschitz-pt1}.
    The base of induction is $\ell=k$
    and holds by the assumption
    of \eqref{prop:ml-lipschitz-pt2}.
    Assume now that $F$ is $k$-increasing on each $\ell$-face for some $\ell$ with  $\ell\leq d-1$ and prove that it is 
    $k$-increasing on each $(\ell+1)$-face.
    Let
    $\cB_{\uu,A}$ be an arbitrary $(\ell+1)$-face of $\cB$, where $\uu\in \vertt \cB$
    and $A\subseteq \{1,2,\ldots,d\}$ with $\Card A=\ell+1$.
    Let $\uz\in \cB_{\uu,A}$ and 
    $(\cB_{\uu,A})_{\uz,B}$ be a $\ell$-section of $\cB_{\uu,A}$, where $\uz\in \cB_{\uu,A}$ and 
    $B\subset A$ with $\Card B=\ell$. Let $\{i_0\}:=A\setminus B$. Without loss of generality we can assume that $\uu_{i_0}=x_{u_0}$. Then 
    $z_{i_0}=\lambda x_{i_0}+(1-\lambda)y_{i_0}$ for some 
    $\lambda\in [0,1]$ and thus the $l$ section $(\cB_{\uu, A})_{\uz, B}$ can be written as a convex combination of
    two parallel $l$ sub-faces $\cB_{\uu, B}$ and $\cB_{\uv, B}$
    \begin{equation}
    \label{multilinear}
    (\cB_{\uu,A})_{\uz,B}=
    \lambda \cB_{\uu,B}
    +
    (1-\lambda) \cB_{\uv,B},
    \end{equation}
    where $\uv$ is a neighbor vertex of $\uu$  that only differs in coordinate $i_0$
    $$
    v_j=\left\{
    \begin{array}{lr}
         u_j&  j\neq i_0,\\
         y_{i_0},&   j=i_0.
    \end{array}
    \right.
    $$ By
    the induction hypothesis, $F$
    is $k$-increasing on $\cB_{\underline u,B}$
    and $\cB_{\underline v,B}$. B by multilinearity 
    of $F$ and by \eqref{multilinear}, $F$ is $k$-increasing on  
    $\cB_{\underline a,A}$.
    This proves the induction step and concludes the proof of the lemma.
\end{proof}

Now we can prove Proposition \ref{prel:k-prop}.

\begin{proof}[Proof of Proposition \ref{prel:k-prop}]The linear program \eqref{LP-general} is obtained by extending the linear program
for general quasi-copulas introduced in \cite[Proposition~2]{OmVuZa25}
with additional constraints enforcing the \(k\)-increasing property.
By \cite[Lemma~2]{AGMB17}, any \(k\)-increasing quasi-copula is also
\(j\)-increasing for every \(j \le k\).
Nevertheless, in the linear program it is necessary to impose the
\(j\)-increasing constraints explicitly for all \(j \le k\).
Otherwise, the resulting quasi-copula would satisfy the
\(k\)-increasing property only on the box
\(\prod_{i=1}^k [a_i,b_i]\), which is insufficient to invoke
\cite[Lemma~2]{AGMB17}.

It remains to verify that the solution of the linear program \eqref{LP-general} can be extended to a $k$-increasing quasi-copula.

Let points $x_\mi\in \prod_{i=1}^d\{a_i,1\}=:\calD$ be defined by \eqref{def:xI-v2}
and real numbers 
$q_\mi$ for all $\mi\in \{0,1\}^d$,
satisfy conditions \eqref{LP-general}. 
Let $\calL_i$ denote the $(d-1)$--dimensional faces of $[0,1]^d$ containing $(0,\ldots,0)$, i.e., 
\begin{equation*}
\calL_i=
\{
(x_1,\ldots,x_{i-1},0,x_{i+1},\ldots,x_d)
\colon x_i\in [0,1]
\}.
\end{equation*}
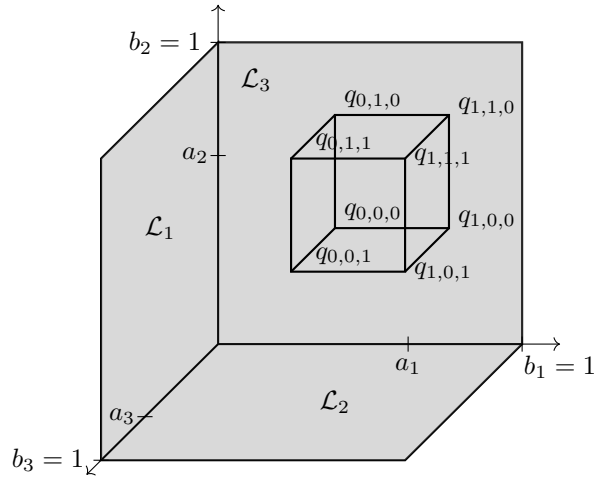
\begin{figure}[h!]
        \centering
    
        \begin{tikzpicture}
        
        \pgfmathsetmacro{\one}{4}
        \pgfmathsetmacro{\aa}{2.5}
        \pgfmathsetmacro{\ab}{2.5}
        \pgfmathsetmacro{\ac}{2.5}
        \pgfmathsetmacro{\ba}{4}
        \pgfmathsetmacro{\bb}{4}
        \pgfmathsetmacro{\bc}{4}

        \usetikzlibrary{patterns}
        \draw[thick] (0, 0, 0) -- (\one, 0, 0) -- (\one, \one, 0) -- (0,\one, 0);
        \draw[thick] (0, 0, 0) -- (0, \one, 0) -- (0, \one, \one) -- (0, 0,\one);
        \draw[thick] (0, 0, 0) -- (0, 0, \one) -- (\one, 0, \one) -- (\one, 0,  0);
        \fill[gray,opacity=0.3] (0, 0, 0) -- (\one, 0, 0) -- (\one, \one, 0) -- (0,\one, 0);
        \fill[gray,opacity=0.3] (0, 0, 0) -- (0, \one, 0) -- (0, \one, \one) -- (0, 0,\one);
        \fill[gray,opacity=0.3] (0, 0, 0) -- (0, 0, \one) -- (\one, 0, \one) -- (\one, 0,  0);
        \node at (0.5, 3.5, 0) {$\calL_3$};
        \node at (2.5, 0.2, 2.5) {$\calL_2$};
        \node at (0.2, 2.5, 2.5) {$\calL_1$};
        
        \coordinate (A) at (\aa, \ab, \ac); 
        \coordinate (B) at (\ba, \ab, \ac); 
        \coordinate (C) at (\ba, \bb, \ac); 
        \coordinate (D) at (\aa, \bb, \ac); 
        \coordinate (E) at (\aa, \ab, \bc); 
        \coordinate (F) at (\ba, \ab, \bc); 
        \coordinate (G) at (\ba, \bb, \bc); 
        \coordinate (H) at (\aa, \bb, \bc); 
        
        \draw[thick] (A) -- (B) -- (C) -- (D) -- cycle;
        \draw[thick] (E) -- (F) -- (G) -- (H) -- cycle;
        \draw[thick] (A) -- (E);
        \draw[thick] (B) -- (F);
        \draw[thick] (C) -- (G);
        \draw[thick] (D) -- (H);
        
        \node at (\aa + 0.5, \ab + 0.2, \ac) {$q_{0,0,0}$};
        \node at (\ba + 0.5, \ab + 0.1, \ac) {$q_{1,0,0}$};
        \node at (\aa + 0.5, \bb + 0.2, \ac) {$q_{0,1,0}$};
        \node at (\ba + 0.5, \bb + 0.1, \ac) {$q_{1,1,0}$};
        \node at (\aa + 0.7, \ab + 0.2, \bc)  {$q_{0,0,1}$};
        \node at (\ba + 0.5, \ab + 0.0, \bc){$q_{1,0,1}$};
        \node at (\aa + 0.7, \bb + 0.2, \bc) {$q_{0,1,1}$};
        \node at (\ba +0.5, \bb, \bc) {$q_{1,1,1}$};
        
        \node at (\aa, -0.3, 0) {$a_1$};
        \draw (\aa, -0.1, 0) -- (\aa, 0.1, 0);
        \node at (\ba+0.5, -0.3, 0) {$b_1=1$};
        \draw (\ba, -0.1, 0) -- (\ba, 0.1, 0);
        \node at (-0.3, \ab, 0) {$a_2$};
        \draw (-0.1, \ab, 0) -- (0.1, \ab, 0);
        \node at (-0.7, \bb, 0) {$b_2=1$};
        \draw (-0.1, \bb, 0) -- (0.1, \bb, 0);
        \node at (-0.3, 0, \ac) {$a_3$};
        \draw (-0.1, 0, \ac) -- (0.1, 0, \ac);
        \node at (-0.7, 0, \bc) {$b_3=1$};
        \draw (-0.1, 0, \bc) -- (0.1, 0, \bc);
        
        \draw[->, thin] (0,0,0) -- (4.5,0,0);
        \draw[->, thin] (0,0,0) -- (0,4.5,0);
        \draw[->, thin] (0,0,0) -- (0,0,4.5);
        
        \end{tikzpicture}
    \caption{To construct a $k$-increasing 3-quasi-copula $Q$ given the values $q_\mi$ at $x_\mi\in \calD:=\prod_{i=1}^3\{a_i,1\}$, we first define it to be 0 on all $2$-dimensional faces $\calL_1$, $\calL_2,$ $\calL_3$ and prove that the extension indeed
    meets the requirements of a $k$-increasing quasi-copula.
    }
    \label{fig:facesL}
\end{figure}

Let
\begin{equation*} 
\cD^{(\ext)}
    :=
    \cD 
    \cup 
    \Big(\cup_{i=1}^d \calL_i\Big) 
\end{equation*}
and define 
$$
    Q:\calD^{(\ext)}\to\RR
$$ 
by
\begin{equation*}
    Q(x_1,\ldots,x_d)=\left\{
    \begin{array}{rl}
    q_\mi,& \text{if }(x_1,\ldots,x_d)=x_\mi \text{ for some }
        \mi\in \{0,1\}^d,\\
    0,& \text{if }(x_1,\ldots,x_d)\in \calL_i \text{ for some }i\in \{1,\ldots,d\}.
    \end{array}
    \right.
\end{equation*}

We subdivide the box $[0,1]^d$
into $2^d$ smaller $d$-boxes 
\begin{equation}
\label{subboxes}
    \cB_{\mi}=\prod_{j=1}^d \delta_j(\mi)
\end{equation}
for $\mi=(\mi_1,\ldots,\mi_d)\in \{0,1\}^d$, where
    $$
    \delta_j(\mi)
    =
    \left\{
    \begin{array}{rl}
    [0,a_j],&   \text{if }\mi_j=0,\\
    \left[a_{j},1\right],& \text{if }\mi_j=1.
    \end{array}
    \right.
    $$
In particular,
\begin{align*}
    \calB_{(0,0\ldots 0)} &= \prod_{k=1}^d[0, a_k],
    \quad
    \calB_{(1, 0\ldots 0)}=[a_1, 1]\times \prod_{k=2}^d [0, a_k],\ldots,\quad
    \calB_{(1,1\ldots 1)}=\prod_{k=1}^d[a_k, 1].
\end{align*}

Note that the $Q$-volume $V_Q(\calB_\mi)$ of each box $\calB_\mi$ is determined by the value of $Q$
on the vertices of $\calB_\mi$.
For each $\mi\in \{0,1\}^d$ we define a constant function
    $$
    \rho_\mi:\cB_\mi\to \RR,\quad
    \rho_\mi:=\frac{V_Q(\cB_\mi)}{\prod_{j=1}^d V(\delta_j(\mi)))},
    $$
where $V([a, b]) = |b - a|$ is the length of the interval $[a, b]$.
Let us define a piecewise constant function 
\begin{equation*}
\rho:[0,1]^d\to \RR,\quad
\rho(\x):=
\left\{
\begin{array}{rl}
\rho_\mi,&  \text{if }\x\in \Int(\cB_\mi) \text{ for some }\mi\in \{0,1\}^d,\\
0,& \text{otherwise},
\end{array}
\right.
\end{equation*}
where $\Int(A)$ stands for the topological interior of the set $A$ in the usual Euclidean topology.
We will prove that a function $Q:[0,1]^d\to \RR$,
defined by 
\begin{equation}
    Q(x_1, x_2, \ldots x_d) = \int_0^{x_1}\int_0^{x_2}\ldots \int_0^{x_d} 
    \rho(x_1, x_2, \ldots x_d) \dd x_1 \dd x_2\cdots \dd x_d,
    \label{eq:integral}  
\end{equation}
is a $k$--increasing quasi-copula satisfying the statement of the proposition.

By the same proof as in \cite[Theorem 2.1]{BeOmVuZa},
it follows that $Q$ is a quasi-copula. It remains to prove that it is $k$-increasing. By Lemma \ref{multilinear}, it suffices to prove that $Q$
is $k$-increasing on each $k$-face of each box $\cB_{\mi}$ for $\mi\in \{0,1\}^d$. By construction, this is true for 
$\mi_{\text{main}}=(1,1,\ldots,1)$. Let us now consider 
$\mi\in \{0,1\}^d\setminus \mi_{\text{main}}$.
If the face has constant $0$ on some component, then $Q$ is constantly equal to 0 on the whole face and $Q$ is $k$-increasing. Otherwise every constant component is equal to $a_j$ or $1$. Further on, we can assume that all non-constant components are $[a_j,1]$ and not $[0,a_j]$, because otherwise we can further project the section on the $j$-th component to $\{a_j\}$, hence obtaining a lower dimensional face of $\cB_{\mi_{\text{main}}}$. However, this reduction give us a $k$-dimensional face of 
$\cB_{\mi_{\text{main}}}$, whence it is $k$-increasing by construction.
\end{proof}



By the following proposition it suffices to consider symmetric solutions
to the linear program \eqref{LP-general}.

\begin{proposition}
\label{k-symmetric-solution}
Assume the notation from Proposition \ref{prel:k-prop}.
There exists an optimal solution to the linear program \eqref{LP-general} 
of the form
\begin{equation}
\label{symmetric-solution-new}
    \big(
    \underbrace{a^\ast,\ldots,a^{\ast}}_{d},
    \underbrace{b^\ast,\ldots,b^{\ast}}_{d},
    q_{\|\mi^{(1)}\|_1}^\ast,\ldots,q_{\|\mi^{(2^d)}\|_1}^\ast
    \big)
\end{equation}
for some $a,b,q_1,\ldots,q_d\in [0,1]$.

    Analogously, replacing $\min$ with $\max$ in \eqref{LP-general} above,
    the same statement holds.
\end{proposition}

\begin{proof}
    The proof is the same as for \cite[Proposition 3]{BeOmVuZa}, but we include it for the sake of completeness.
    Let $S_d$ be the set of all permutations of a $d$-element set
    $\{1,\ldots,d\}$.
    For $\Phi\in S_d$ and $\mi^{(j)}:=(\mi_1^{(j)},\ldots,\mi_d^{(j)})\in\{0,1\}^d$,
    let $\Phi(\mi^{(j)}):=\big(\mi_{\Phi(1)}^{(j)},\ldots,\mi_{\Phi(d)}^{(j)}\big)$.
    For every optimal solution     
    $\big(
    a_1^\ast,\ldots,a_d^{\ast},
    b_1^\ast,\ldots,b_d^{\ast},
    q_{\mi^{(1)}}^\ast,\ldots,q_{\mi^{(2^d)}}^\ast
    \big)$
    to \eqref{LP-general},
    $\big(
    a_1^\ast,\ldots,a_d^{\ast},
    b_1^\ast,\ldots,b_d^{\ast},
    q_{\mi^{(1)}}^\ast,\ldots,q_{\mi^{(2^d)}}^\ast
    \big)$
    is also an optimal solution to \eqref{LP-general},
    whence an optimal solution  to \eqref{LP-general} of the form 
    \eqref{symmetric-solution-new}
    is equal to 
    $\frac{1}{d!}
    \sum_{\Phi\in S_d}
    \big(
    a_{\Phi(1)}^\ast,\ldots,a_{\Phi(d)}^{\ast},
    b_{\Phi(1)}^\ast,\ldots,b_{\Phi(d)}^{\ast},
    q_{\Phi(\mi^{(1)})}^\ast,\ldots,q_{\Phi(\mi^{(2^d)})}^\ast
    \big).
    $
\end{proof}

An immediate corollary to Proposition \ref{k-symmetric-solution} is the following.

\begin{corollary}
\label{cor:sim-k}
    The optimal value of the linear program \eqref{LP-general} is equal to the optimal value of the 
    linear program
\begin{align}
\label{LP-symmetric-k}
\begin{split}
\min_{
\substack{
    a,
    b,
    q_0,q_{1},\ldots,q_d
}}
&\hspace{0.5cm} \sum_{i=0}^d (-1)^{d-i} \binom{d}{i} q_i,\\
\text{subject to }
&\hspace{0.5cm} 
    0\leq a< b\leq 1,\\
&\hspace{0.5cm} 
     q_{i} - q_{i-1} \le b - a\quad \text{for }i=1,2,\ldots,d,\\
&\hspace{0.5cm} 
    0\leq \sum_{i=0}^j (-1)^{j-i}\binom{j}{i}q_{\ell-i}
    \quad
    \text{for all }j\in 1,2,\ldots,k\text{ and all }\ell=j,j+1,\ldots,d,\\
&\hspace{0.5cm} 
    \max\{0,(d-i)a+ib-d+1\} 
    \le q_i \le 
    a\quad \text{for } i= 0,1,\ldots,d-1,\\
&\hspace{0.5cm} 
    \max\{0,db-d+1\} 
    \le q_d \le 
    b.    
\end{split}
\end{align}

    Analogously, replacing $\min$ with $\max$ in \eqref{LP-symmetric-k} above,
    the same statement holds.
\end{corollary}

It turns out that some of the constraints in \eqref{LP-symmetric-k} are redundant, while introducing new variables 
\begin{align}
\label{new-variables-level-k}
\begin{split}
    \delta_i^{(1)}&:=q_i-q_{i-1}\quad \text{for }i=1,2,\ldots,d,\\
    \delta_i^{(j)}&:=\delta_i^{(j-1)}-\delta_{i-1}^{(j-1)}
        \quad\text{for } j\in 2,3,\ldots,k\text{ and }i=j,j+1,\ldots,d.
\end{split}
\end{align}
further decreases the number of constraints.
Let us also define the numbers 
    $\beta_{i}^{(j)}$ 
recursively by
\begin{align}
\label{recursiveness}
\begin{split}
    \beta_{i}^{(j)}&=
        \left\{
        \begin{array}{rl}
            d-i & \text{for }j=2\text{ and }i= 2,3,\ldots,d-1,\\[0.2em]
            \beta_i^{(j-1)}& \text{for }j=3,4,\ldots,k\text{ and }i=2,3,\ldots,j-2,\\[0.2em]
            \sum_{\ell=i}^{d-1}\beta_{\ell}^{(j-1)}& 
                \text{for }j=3,4,\ldots,k\text{ and }i=j-1,j,\ldots,d-1,
        \end{array}
        \right..
\end{split}
\end{align}

\begin{lemma}
\label{alpha+beta=gamma}
    Let $\alpha_{i}^{(j)}$, $\beta_i^{(j)}$ and $\gamma_i^{(j)}$
    be as in \eqref{recursiveness-v2},
    \eqref{recursiveness-v3} and \eqref{recursiveness}.
    For $j= 2,3,\ldots,k$ and $i= 2,3,\ldots,d$, we have
        \begin{equation}
            \label{eq:alpha+beta=gamma}
                \alpha_{i}^{(j)}+\beta_{i}^{(j)}=\gamma_{i}^{(j)}.
        \end{equation}
\end{lemma}

\begin{proof}
    For $j=2$, we have 
        $\alpha_i^{(2)}+\beta_i^{(2)}=d+1-i=\gamma_i^{(2)}$ for $i= 2,3,\ldots,d$,
    which proves \eqref{eq:alpha+beta=gamma}. 
    For $j>2$, \eqref{eq:alpha+beta=gamma} follows inductively using that 
    \eqref{eq:alpha+beta=gamma} holds for smaller values of $j$.
\end{proof}

\begin{proposition}
\label{simplify-further-general}
    Let $\alpha_{i}^{(j)}$, $\beta_i^{(j)}$ and $\gamma_i^{(j)}$
    be as in \eqref{recursiveness-v2},
    \eqref{recursiveness-v3}
    and \eqref{recursiveness}.
    The optimal value of the linear program \eqref{LP-symmetric-k} is equal to the optimal value of the linear program
\begin{align}
\label{LP-symmetric-general}
\begin{split}
\min_{
\substack{
    a,
    b,
    q_0,\\
    \delta_{1}^{(1)},\delta_{2}^{(2)},\ldots,\delta_{k-1}^{(k-1)},\\
    \delta_k^{(k)},\delta_{k+1}^{(k)},\ldots,\delta_d^{(k)}
}}
&\hspace{0.5cm}
\sum_{j=k}^d (-1)^{d+j}\binom{d-k}{j-k}\delta^{(k)}_j,\\
\text{subject to }
&\hspace{0.5cm} 
    b\leq 1,\\
&\hspace{0.5cm} 
    a-b+\delta_1^{(1)}+
        \sum_{i=2}^{k-1}\alpha_{i}^{(k)}\delta^{(i)}_i+
        \sum_{i=k}^{d}\alpha_{i}^{(k)}\delta^{(k)}_i\le 0,\\
&\hspace{0.5cm} 
    -a+q_0+(d-1)\delta^{(1)}_1+
    \sum_{i=2}^{k-1}\beta_{i}^{(k)}\delta^{(i)}_i+
        \sum_{i=k}^{d-1}\beta_{i}^{(k)}\delta^{(k)}_i
    \leq 0,\\
&\hspace{0.5cm} 
    db-q_0-d\delta_1^{(1)}
    -\sum_{i=2}^{k-1}\gamma_{i}^{(k)}\delta^{(i)}_i
    -\sum_{i=k}^{d}\gamma_{i}^{(k)}\delta^{(k)}_i
    \leq d-1,\\
&\hspace{0.5cm} 
    a\geq 0,\; b\geq 0,\; q_0\geq 0,\; 
    \delta_i^{(i)}\geq 0\;\; \text{for }i=1,2,\ldots,k-1,\\
&\hspace{0.5cm}     
    \delta_i^{(k)}\geq 0\;\; \text{for }i=k,k+1,\ldots,d.
\end{split}
\end{align}

    Analogously, replacing $\min$ with $\max$ in \eqref{LP-symmetric-general} above,
    the same statement holds.
\end{proposition}

\begin{proof}
We will use induction on $k$ to prove Proposition \ref{simplify-further-general}.

We start by $k=2$.
By \cite[Proposition 5]{OmVuZa25},  \eqref{LP-symmetric-k}
is equivalent to
\begin{align}
\label{LP-symmetric-level-k-v2}
\begin{split}
\min_{
\substack{
    a,
    b,
    q_0,\delta_{1}^{(1)},\ldots,\delta_d^{(1)}
}}
&\hspace{0.5cm}
\sum_{j=1}^d (-1)^{d+j}\binom{d-1}{j-1}\delta_j^{(1)},\\
\text{subject to }
&\hspace{0.5cm} 
    b\leq 1,\\
&\hspace{0.5cm} 
     \delta_i^{(1)} \le b - a\quad \text{for }i=1,2,\ldots,d,\\
&\hspace{0.5cm} 
     0\le q_{i} - 2q_{i-1}+q_{i-2} \quad \text{for }i=2,3,\ldots,d,\\
&\hspace{0.5cm} 
    q_0+\sum_{i=1}^{d-1}\delta_i^{(1)}\leq a\\
&\hspace{0.5cm} 
    db-d+1 
    \le q_0+\sum_{i=1}^{d}\delta_i^{(1)}\\
&\hspace{0.5cm} 
    a\geq 0,\; b\geq 0,\; q_0\geq 0,\; \delta_i^{(1)}\geq 0 \quad \text{for }i=1,2,\ldots,d.  
\end{split}
\end{align}
Since $\delta_i^{(1)}, \delta_{i}^{(2)}$ are as in \eqref{new-variables-level-k},
we can replace the inequalities 
    $$0\le q_{i} - 2q_{i-1}+q_{i-2} \quad \text{for }i=2,3,\ldots,d$$
with the inequalities
\begin{equation} 
    \label{new-inequalities}
        0\leq \delta_i^{(2)}\quad \text{for }i=2,3,\ldots,d.
\end{equation}
Using \eqref{new-inequalities},
the constraints
$\delta_i^{(1)}\geq 0$ and $\delta_i^{(1)}\leq b-a$ for $i=1,2,\ldots,d$ in \eqref{LP-symmetric-level-k-v2}
can be replaced 
by the constraints
$\delta_1^{(1)}\geq 0$ and $\delta_d^{(1)}\leq b-a.$
Since
    $\delta_j^{(1)}=\delta_j^{(2)}+\delta_{j-1}^{(1)}$
for $j=2,3,\ldots,d$, it follows inductively that
\begin{equation}
\label{061025-2113}
    \delta_j^{(1)}=\sum_{l=2}^{j}\delta_l^{(2)}+\delta_{1}^{(1)}
    \quad \text{for } j=2,3,\ldots,d.
\end{equation}
Using \eqref{061025-2113} in the constraints of \eqref{LP-symmetric-level-k-v2}, it is easy to see that we get the constraints of \eqref{LP-symmetric-k}
with $\alpha_i^{(2)},\beta_{i}^{(2)}, \gamma_i^{(2)}$. 
It remains to show that the objective function of \eqref{LP-symmetric-k}
becomes the objective function of \eqref{LP-symmetric-k} after substitutions \eqref{061025-2113}. Indeed,
\begin{align*}
\sum_{j=1}^d (-1)^{d-j}\binom{d-1}{j-1}\delta_{j}^{(1)}
&=
\sum_{j=1}^d (-1)^{d-j}\binom{d-1}{j-1}
\Big(
\delta_1^{(1)}+\sum_{\ell=2}^{j}\delta_{\ell}^{(2)}
\Big)\\
&=
\delta_{1}^{(1)}\underbrace{\sum_{j=1}^d (-1)^{d-j}\binom{d-1}{j-1}}_{0}
+
\sum_{\ell=2}^d \delta_{\ell}^{(2)}
    \underbrace{\Big(
        \sum_{k=\ell}^d (-1)^{d-k}\binom{d-1}{k-1}
    \Big)}_{
    (-1)^{d+\ell}\binom{d-2}{\ell-2}
    },
\end{align*}
where in the last equality we used that for $0\leq n\leq r$ we 
have
\begin{equation} 
    \label{alternating-binomial}
        \sum_{k= n}^{r}(-1)^k\binom{r}{k}
        \underbrace{=}_{\sum_{k= 0}^{r}(-1)^k\binom{r}{k}=0}-\sum_{k=0}^{n-1}(-1)^k\binom{r}{k}=(-1)^n\binom{r-1}{n-1}.
\end{equation}

Assume know that the proposition holds for some $k$, $2\leq k\leq d-1$ and prove it for $k+1$. By the induction hypothesis the optimal value of the linear program \eqref{LP-symmetric-k} is equal to the optimal value of the linear program
\begin{align}
\label{LP-symmetric-general-v2}
\begin{split}
\min_{
\substack{
    a,
    b,
    q_0,\\
    \delta_{1}^{(1)},\delta_{2}^{(2)},\ldots,\delta_{k-1}^{(k-1)},\\
    \delta_k^{(k)},\delta_{k+1}^{(k)},\ldots,\delta_d^{(k)}
}}
&\hspace{0.5cm}
\sum_{j=k}^d (-1)^{d+j}\binom{d-k}{j-k}\delta^{(k)}_j,\\
\text{subject to }
&\hspace{0.5cm} 
    b\leq 1,\\
&\hspace{0.5cm} 
    a-b+\delta_1^{(1)}+
        \sum_{i=2}^{k-1}\alpha_{i}^{(k)}\delta^{(i)}_i+
        \sum_{i=k}^{d}\alpha_{i}^{(k)}\delta^{(k)}_i\le 0,\\
&\hspace{0.5cm} 
    0\leq \sum_{i=0}^{k+1} (-1)^{k+1-i}\binom{k+1}{i}q_{\ell-i}
    \quad
    \text{for all }\ell=k+1,k+2,\ldots,d,\\
&\hspace{0.5cm} 
    -a+q_0+(d-1)\delta^{(1)}_1+
    \sum_{i=2}^{k-1}\beta_{i}^{(k)}\delta^{(i)}_i+
        \sum_{i=k}^{d-1}\beta_{i}^{(k)}\delta^{(k)}_i
    \leq 0,\\
&\hspace{0.5cm} 
    db-q_0-d\delta_1^{(1)}-\sum_{i=2}^{k-1}\gamma_{i}^{(k)}\delta^{(i)}_i-
        \sum_{i=k}^{d}\gamma_{i}^{(k)}\delta^{(k)}_i
    \leq d-1,\\
&\hspace{0.5cm} 
    a\geq 0,\; b\geq 0,\; q_0\geq 0,\; 
    \delta_i^{(i)}\geq 0\;\; \text{for }i=1,2,\ldots,k-1,\\
&\hspace{0.5cm} 
    \delta_i^{(k)}\geq 0\;\; \text{for }i=k,k+1,\ldots,d.
\end{split}
\end{align}
Since $\delta_i^{(j)}$ are as in \eqref{new-variables-level-k},
we can replace the inequalities 
    $$0\leq \sum_{i=0}^{k+1} (-1)^{k+1-i}\binom{k+1}{i}q_{\ell-i}
    \quad
    \text{for all }\ell=k+1,k+2,\ldots,d,$$
with the inequalities
    $$0\leq \delta_\ell^{(k+1)}\quad \text{for }
    \ell= k+1,k+2,\ldots,d.$$
Since
    $\delta_j^{(k)}=\delta_j^{(k+1)}+\delta_{j-1}^{(k)}$
for $j=k+1,k+2,\ldots,d$, it follows inductively that
\begin{equation}
\label{171025-1847}
    \delta_j^{(k)}=\sum_{\ell=k+1}^{j}\delta_\ell^{(k+1)}+\delta_{k}^{(k)}
    \quad \text{for } j=k+1,k+2,\ldots,d.
\end{equation}
Using \eqref{171025-1847} in the constraints of \eqref{LP-symmetric-general-v2}, it is easy to see that the constraints become as in \eqref{LP-symmetric-general}, where $k$ is replaced by $k+1$.
It remains to show that the objective function of \eqref{LP-symmetric-general-v2}
becomes the objective function of \eqref{LP-symmetric-general} after substitutions \eqref{171025-1847}. Indeed,
\begin{align*}
\sum_{j=k}^d (-1)^{d+j}\binom{d-k}{j-k}\delta^{(k)}_j
&=
\sum_{j=k}^d (-1)^{d+j}\binom{d-k}{j-k}
\Big(
\delta_k^{(k)}+\sum_{\ell=k+1}^{j}\delta_{\ell}^{(k+1)}
\Big)\\
&=
\delta_{k}^{(k)}\underbrace{\sum_{j=k}^d (-1)^{d+j}\binom{d-k}{j-k}}_{0}
+
\sum_{\ell=k+1}^d \delta_{\ell}^{(k+1)}
    \underbrace{\Big(
        \sum_{j=\ell}^d (-1)^{d+j}\binom{d-k}{j-k}
    \Big)}_{
    (-1)^{d+\ell}\binom{d-k-1}{\ell-k-1}
    },
\end{align*}
where in the last equality we used \eqref{alternating-binomial}.
\end{proof}

Using duality, the following proposition holds.

\begin{proposition}
\label{simplify-further-general-v2}
    Let $\alpha_{i}^{(j)}$, $\beta_i^{(j)}$ and $\gamma_i^{(j)}$
    be as in \eqref{recursiveness-v2},
    \eqref{recursiveness-v3}
    and \eqref{recursiveness}.
    The optimal value of the linear program \eqref{LP-symmetric-general} is equal to the optimal value of the linear program
\begin{align}
\label{dual-LP-symmetric-general}
\begin{split}
\max_{
\substack{
    y_1,y_2,y_3,y_4
}}
&\hspace{0.5cm} -y_1-(d-1)y_4,\\
\text{subject to }
&\hspace{0.5cm} 
    y_1-y_2+dy_4\ge 0,\\
&\hspace{0.5cm}
    y_2-y_3\ge 0,\\
&\hspace{0.5cm}
    y_3-y_4\ge 0,\\ 
&\hspace{0.5cm}
    y_2+(d-1)y_3-dy_4\geq 0,\\
&\hspace{0.5cm}
    \alpha_i^{(k)}y_2+\beta_{i}^{(k)}y_3-\gamma_i^{(k)}y_4\geq 0
    \quad\text{for }i=2,3,\ldots,k-1,
    \\
&\hspace{0.5cm}
    \alpha_i^{(k)}y_2+\beta_{i}^{(k)}y_3-\gamma_i^{(k)}y_4\geq (-1)^{d+1+i}\binom{d-k}{i-k}
    \quad\text{for }i=k,k+1,\ldots,d,\\
&\hspace{0.5cm} 
    y_1\geq 0,\; y_2\geq 0,\; y_3\geq 0,\; y_4\geq 0.
\end{split}
\end{align}

Analogously, the optimal value of the linear program \eqref{LP-symmetric-general} with the objective
function
$$\max_{
\substack{
    a,
    b,
    q_0,\\
    \delta_{1}^{(1)},\delta_{2}^{(2)},\ldots,\delta_{k-1}^{(k-1)},\\
    \delta_k^{(k)},\delta_{k+1}^{(k)},\ldots,\delta_d^{(k)}
}}
\sum_{j=k}^d (-1)^{d+j}\binom{d-k}{j-k}\delta^{(k)}_j,
$$
is equal to the optimal value of the linear program
\begin{align}
\label{dual-LP-symmetric-general-v2}
\begin{split}
\min_{
\substack{
    y_1,y_2,y_3,y_4
}}
&\hspace{0.5cm} y_1+(d-1)y_4,\\
\text{subject to }
&\hspace{0.5cm} 
    y_1-y_2+dy_4\ge 0,\\
&\hspace{0.5cm}
    y_2-y_3\ge 0,\\
&\hspace{0.5cm}
    y_3-y_4\ge 0,\\ 
&\hspace{0.5cm}
    y_2+(d-1)y_3-dy_4\geq 0,\\
&\hspace{0.5cm}
    \alpha_i^{(k)}y_2+\beta_{i}^{(k)}y_3-\gamma_i^{(k)}y_4\geq 0
    \quad\text{for }i=2,3,\ldots,k-1,
    \\
&\hspace{0.5cm}
    \alpha_i^{(k)}y_2+\beta_{i}^{(k)}y_3-\gamma_i^{(k)}y_4\geq (-1)^{d+i}\binom{d-k}{i-k}
    \quad\text{for }i=k,k+1,\ldots,d,\\
&\hspace{0.5cm} 
    y_1\geq 0,\; y_2\geq 0,\; y_3\geq 0,\; y_4\geq 0.
\end{split}
\end{align}
\end{proposition}

\begin{proof}
Note that \eqref{dual-LP-symmetric-general} is the dual linear program to \eqref{LP-symmetric-general}, whence the statement of Proposition \ref{simplify-further-general-v2} follows by the strong duality.
\end{proof}

\section{Proofs of Theorems \ref{sol:general-increasing} and
\ref{sol:general-increasing-maximal}}
\label{sec:proofs}

Finally we can prove Theorem \ref{sol:general-increasing}.

\begin{proof}[Proof of Theorem \ref{sol:general-increasing}]
By the results above to prove Theorem \ref{sol:general-increasing}, it suffices to find the optimal value of
\eqref{dual-LP-symmetric-general}.
First we prove the following claim.\\

\noindent \textbf{Claim 1.} Their exists an optimal solution 
 $(y^{\opt}_1,y^{\opt}_2,y^{\opt}_3,y^{\opt}_4)$ 
to the linear program 
    \eqref{dual-LP-symmetric-general}
which
satisfies
    $y_2^{\opt}=y_3^{\opt}.$\\

\noindent \textit{Proof of Claim 1.}
Let
    $(y^{\opt}_1,y^{\opt}_2,y^{\opt}_3,y^{\opt}_4)$ 
be an optimal solution to \eqref{dual-LP-symmetric-general}.
If $y^{\opt}_2>y^{\opt}_3$, then replacing $y^{\opt}_3$ with $y^{\opt}_2$
it is clear that $(y^{\opt}_1,y^{\opt}_2,y^{\opt}_2,y^{\opt}_4)$ satisfies all constraints and the value of the objective function remains unchanged.\hfill$\blacksquare$\\

Using Claim 1 and \eqref{eq:alpha+beta=gamma},
the linear program \eqref{dual-LP-symmetric-general} simplifies to:

\begin{align}
\label{dual-LP-symmetric-general-v2}
\begin{split}
\max_{
\substack{
    y_2,y_4
}}
&\hspace{0.5cm} -y_1-(d-1)y_4,\\
\text{subject to }
&\hspace{0.5cm} 
    y_1-y_2+dy_4\ge 0,\\
&\hspace{0.5cm}
    y_2-y_4\ge 0,\\ 
&\hspace{0.5cm}
    \gamma_i^{(k)}(y_2-y_4)\geq (-1)^{d+1-i}\binom{d-k}{i-k}\quad \text{for }i=k,k+1,\ldots,d,\\
&\hspace{0.5cm} 
    y_1\geq 0,\; y_2\geq 0,\; y_4\geq 0.
\end{split}
\end{align}

\noindent \textbf{Claim 2.} There is an optimal solution 
$(y^{\opt}_1,y^{\opt}_2,y^{\opt}_4)$ 
to the linear program 
    \eqref{dual-LP-symmetric-general-v2} 
such that 
    $$y_2^\ast=y_1^\ast+dy_4^\ast.$$

\noindent \textit{Proof of Claim 2.}
Let
    $(y^{\ast}_1,y^{\opt}_2,y^{\opt}_4)$ 
be an optimal solution to \eqref{dual-LP-symmetric-general-v2}.
If $y_2^\ast<y_1^\ast+dy_4^\ast$, then we can replace $y_2^\ast$ with $y_1^\ast+dy_4^\ast$, still satisfying all constraints of
\eqref{dual-LP-symmetric-general-v2} and not changing its objective function. This proves Claim 2.
\hfill$\blacksquare$\\

Using Claim 2, \eqref{dual-LP-symmetric-general-v2} simplifies to

\begin{align}
\label{dual-LP-symmetric-general-v3}
\begin{split}
\max_{
\substack{
    y_1,y_4
}}
&\hspace{0.5cm} -y_1-(d-1)y_4,\\
\text{subject to }
&\hspace{0.5cm}
    \gamma_i^{(k)}(y_1+(d-1)y_4)\geq (-1)^{d+1-i}\binom{d-k}{i-k}\quad \text{for }i= k,k+1,\ldots,d,\\
&\hspace{0.5cm} 
    y_1\geq 0,\; y_4\geq 0.
\end{split}
\end{align}

Clearly, the optimal solution satisfies
\begin{equation}
\label{optimal-sol-general}
y_1^\ast+(d-1)y_4^\ast=
        \max\big(0,\max_{i\in k,\ldots,d}
        (-1)^{d+1-i}\binom{d-k}{i-k}(\gamma_i^{(k)})^{-1}\big),
\end{equation}
which proves the first part of Theorem \ref{sol:general-increasing}.

It remains to prove the part about the realization of one optimal solution. Define
$$q_0^\ast=
(\delta_1^{(1)})^{\ast}=
(\delta_2^{(2)})^{\ast}=
\ldots=
(\delta_{k-1}^{(k-1)})^{\ast}=
(\delta_k^{(k)})^{\ast}=
\ldots=
(\delta_{i_0-1}^{(k)})^{\ast}=
(\delta_{i_0+1}^{(k)})^{\ast}=
\ldots=
(\delta_{d}^{(k)})^{\ast}=0.
$$
Since there is a solution to \eqref{dual-LP-symmetric-general}
such that $y_i^{\ast}\neq 0$ for $i=1,2,3,4$ (This is due to the fact that any pair $(y_1,y_4)\in \RR_{\geq 0}\times \RR_{\geq 0}$ such that 
$y_1+(d-1)y_4=w_{k,d,-}$, is an optimal solution.),
it follows by complementary slackness that in the optimal solution to \eqref{LP-symmetric-general} we have four equalities.
From the optimal value of the objective function of \eqref{LP-symmetric-general} we conclude 
$(\delta_{i_0}^{(k)})^\ast=(\gamma_{i_0}^{(k)})^{-1}$. From the four equalities stated above we conclude \eqref{one:realiziation} using also that
$\delta_j^{(i)}=\sum_{\ell=0}^{i}\binom{i}{\ell}(-1)^\ell q_{j-\ell}.$
\end{proof}

\begin{proof}[Proof of Theorem \ref{sol:general-increasing-maximal}]
The proof is analogous to the proof of Theorem \ref{sol:general-increasing}.
\end{proof}

\section{Concluding remarks and future research}
\label{sec:concluding}
In this paper, we addressed the problem of determining the extreme values of the mass distribution associated with a multidimensional $k$-increasing $d$-variate quasi-copula for $k = 2,3,\ldots,d$, building on the linear programming approach that solved the case $k = 1$ (see \cite{BMUF07} for $d = 3$, \cite{UF23} for $d = 4$, \cite{BeOmVuZa} for $d \leq 17$, and \cite{OmVuZa25} for the general case). We derived closed formulas for the extreme volumes as well as an explicit realization of one corresponding quasi-copula. Besides the symmetrization trick from \cite{OmVuZa25}, the main novelty enabling us to solve the linear programs was the introduction of new sets of variables, which significantly reduced their complexity. We also provided a numerical analysis of the results for dimensions $d \leq 15$.

Finally, we outline several directions for future research.
In the bivariate case, the 2-increasing property coincides with the supermodularity of a function $F$. Very recently, Anzilli and Durante \cite{AD26} derived bounds on the average $F$-volume of closed rectangles in $[0,1]^2$ for a supermodular aggregation function $F$. In line with these results, it would be interesting to study the behavior of extreme volumes for supermodular quasi-copulas or supermodular aggregation functions \cite{DSPS08} in higher dimensions, as well as for ultramodular and modular quasi-copulas investigated in \cite{KMM11}.

It is clear that every extreme point (in the Krein--Milman sense) of the convex set of $k$-increasing quasi-copulas must be a maximal-volume one. Our construction of the extremal maximal-volume solution is symmetric---due to the symmetrization step---and is therefore not extreme. To better understand the geometry of the convex set of $k$-increasing quasi-copulas, it would be natural to characterize their extreme points. For semilinear copulas, such a characterization was obtained in \cite{DFSUF20}.

 \end{document}